\documentclass
{amsart}

\usepackage{amsmath, amsthm, amsfonts, amssymb,mathrsfs}
\usepackage{color,umoline}
\usepackage{enumitem}
\usepackage[all]{xy}
\usepackage{graphicx}
\usepackage{pb-diagram}
\usepackage{caption}
\usepackage{subcaption}
\usepackage{afterpage} 
\usepackage{tikz}
\usepackage{amscd}
\usepackage{txfonts}
\usepackage{relsize}

\usepackage[hidelinks]{hyperref}

\hypersetup{
colorlinks=true,       
linkcolor=blue,          
citecolor=red,        
filecolor=magenta,      
urlcolor=cyan           
}

\newcommand{\Be}{\begin{equation}}
\newcommand{\Ee}{\end{equation}}
\newcommand{\Bea}{\begin{eqnarray}}
\newcommand{\Eea}{\end{eqnarray}}
\newcommand{\Bel}{\begin{align}}
\newcommand{\Eel}{\end{align}}
\newcommand{\Beas}{\begin{eqnarray*}}
\newcommand{\Eeas}{\end{eqnarray*}}
\newcommand{\Benu}{\begin{enumerate}}
\newcommand{\Eenu}{\end{enumerate}}
\newcommand{\Bi}{\begin{itemize}}
\newcommand{\Ei}{\end{itemize}}

\def\norm#1{\|#1\|}
\def\normb#1{\Big\|#1\Big\|}

\def\normo#1{\Big\|#1\Big\|}
\def\Normo#1{\left\|#1\right\|}
\def\wt#1{\widetilde{#1}}
\def\wh#1{\widehat{#1}}

\newcommand{\B}{\mathcal{B}}

\newcommand{\F}{\mathcal{F}}

\newcommand{\C}{\mathbb{C}}

\newcommand{\R}{\mathbb{R}}
\newcommand{\Z}{\mathbb{Z}}

\newcommand{\la}{\lambda}

\newcommand{\ta}{\tau}

\newcommand{\im}{\mathop{\mathrm{Im}}}

\newcommand{\supp}{\operatorname{supp}}

\newcommand{\Del}[1]{}

\newcommand{\finv}[1]{{\mathlarger{ \mathcal F}}^{\,-1}\big({#1}\big)}

\newcommand{\mpq}[1]{\mathfrak M_{p, q} [#1] }
\newcommand{\mpqq}[3]{\mathfrak M_{#1, #2} [#3]}
\newcommand{\Mpqq}[3]{\mathfrak M_{#1, #2} \Big[#3\Big] }

\newcommand{\fQ}{\mathfrak Q}
\newcommand{\fG}{\mathfrak G}
\newcommand{\fS}{\mathfrak S}

\newcommand{\ovpsi}{\overline{\psi}}
\renewcommand{\B}{\Big}
\newcommand{\wave}{\mathlarger{\mathlarger{{\square}}}}
\newcommand{\dn}{\mathbb D}

\numberwithin{equation}{section}

\newtheorem{thm}{Theorem}[section]
\newtheorem{cor}[thm]{Corollary}
\newtheorem{lem}[thm]{Lemma}
\newtheorem{prop}[thm]{Proposition}

\theoremstyle{remark}

\begin{document}

\title[Carleman estimate and Multiplier operator]{Carleman estimates and boundedness of \\ associated multiplier operators}

\author[E. Jeong]{Eunhee Jeong}
\author[Y. Kwon]{Yehyun Kwon}
\author[S. Lee]{Sanghyuk Lee}

\address{Department of Mathematical Sciences, Seoul National University, Seoul 151-747, Republic of Korea}
\email{moonshine10@snu.ac.kr}
\email{kwonyh27@snu.ac.kr}
\email{shklee@snu.ac.kr}

\thanks{This work was supported by NRF of Korea (grant no. 2015R1A2A2A05000956). 
}

\subjclass[2010]{42B15, 35B60} 
\keywords{Carleman estimate, unique continuation}

\begin{abstract}  Let $P(D)$ be the Laplacian $\Delta,$ or the wave operator  $\wave$.    The following type of Carleman estimate is known to be true on a certain range of $p,q$: 
\[\norm{e^{v\cdot x}u}_{L^q(\R^d)} \le C\norm{e^{v\cdot x}P(D)u}_{L^p(\R^d)} \]
with $C$ independent of $v\in \R^d$. The estimates are  consequences of the uniform Sobolev type estimates for second order differential operators  due to Kenig-Ruiz-Sogge \cite{KRS} and Jeong-Kwon-Lee \cite{JKL}.  The  range of $p,q$ for which  the uniform Sobolev type estimates hold was completely characterized for the second order differential operators with nondegenerate principal part. But the optimal range of $p,q$ for which the Carleman estimate holds has not been clarified before.   When $P(D)=\Delta$, $\wave$, or  the heat operator,   we obtain a complete characterization of  the admissible $p,q$  for the aforementioned  type of Carleman estimate.  For this purpose  we investigate $L^p$--$L^q$ boundedness of related multiplier operators. As applications, we also obtain some unique continuation results. 
\end{abstract}
\maketitle

\section{Introduction and results}
In this note we consider Carleman estimates for second order differential operators with a special type of exponential weight. Firstly, we are concerned with  the following  type of Carleman inequality for the Laplacian: 
\begin{equation}\label{car}
\norm{e^{v\cdot x}u}_{L^q(\R^d)} \le C\norm{e^{v\cdot x}(-\Delta)u}_{L^p(\R^d)}
\end{equation}
which holds for all $v\in\R^d$ and $u\in C_0^\infty (\R^d)$ with $C$ independent of $v$. Though, compared with other types of Carleman estimates of nonlinear exponential weights,   the estimate \eqref{car} is relatively simpler to obtain,  {it (or its variants) has various applications.} Especially,  the inequality has been used to study unique continuation properties of differential inequalities, see \cite{KRS, wolff1992, wolff1993}. In particular, it played  an important role  in proving the unique continuation property for the differential inequality 
$|\Delta u|\le V|u|+W|\nabla u|$ (see \cite[Section 6]{wolff1992}).     For $p,q$ on a certain range the estimate \eqref{car} is a consequence of the uniform Sobolev estimate which is due to Kenig, Ruiz, and Sogge \cite{KRS}. In fact,   they  proved that the estimate 
\begin{equation}\label{uni}
\norm{u}_{L^q(\R^d)} \le C\Normo{\left(-\Delta+a\cdot \nabla + b\right)u}_{L^p(\R^d)}
\end{equation}
holds with a constant $C$ independent of  $(a,b)\in\C^d\times\C$ and $u\in W^{2,p}(\R^d)$  if and only if $p,q$ satisfy  
\begin{equation}\label{ran}
\frac 1p-\frac 1q=\frac 2d\,, \quad \frac{d+1}{2d}< \frac 1p < \frac{d+3}{2d}.
\end{equation}
Since $e^{v\cdot x}(-\Delta)e^{-v\cdot x}= -\Delta +2v\cdot\nabla -|v|^2$, for $p,q$ satisfying \eqref{ran} the estimate \eqref{car}  follows from \eqref{uni} by replacing $u$ with $e^{-v\cdot x} u$.

The uniform estimate \eqref{uni} was obtained by using the (seemingly weaker) uniform resolvent estimate  $\| (-\Delta-z)^{-1} u\|_{q}\le C\|u\|_p$, $z\in \C\setminus [0,\infty)$, which is also known to be true if and only if $p,q$ satisfy the same condition \eqref{ran}.  This estimate is closely related to the  Fourier restriction estimate to the sphere. 
A  simple limiting argument (for example, see \cite{KRS, JKL}) shows that the uniform resolvent  estimate implies  the following estimate for the restriction-extension operator defined by the sphere
\begin{equation}\label{rext}
\Normo{\int_{\mathbb S^{d-1}} \wh f (\theta) e^{ix\cdot\theta} d\theta}_{L^q(\R^d)} \lesssim \|f\|_{L^p(\R^d)}.
\end{equation}
This estimate can also be regarded as  an estimate for the Bochner-Riesz operator of order $-1$\footnote{For $\alpha\le -1$, the Bochner-Riesz operator of  order $\alpha$ is defined by $\mathcal F^{-1}(\frac{(1-|\xi|^2)_+^\alpha}{\Gamma(\alpha+1)} \widehat f\,)$ through analytic continuation along $\alpha$.}  (see  \cite{bak, Bo, guti}), and it is well-known that the estimate \eqref{rext} is true only when $ p<2d/(d+1)$ and $q>2d/(d-1)$ (see Theorem \ref{res-ext}).  Combining this with the  condition  $1/p-1/q=2/d$ which is necessary  for \eqref{uni}\footnote{This is easy to see by using the Littlewood-Paley decompositon and Mihlin's multiplier theorem.}, one can see that  the condition \eqref{ran} is also necessary for \eqref{uni} to hold. 


At this point we are naturally led  to question  \emph{whether the range of $p,q$ given in \eqref{ran} is {also} optimal for the Carleman estimate \eqref{car} or there is any other pair of $(p,q)$ for which \eqref{car} is still true}. Clearly, such a question can not be handled by considering the uniform Sobolev inequalities, and 
 it is unlikely  that the range \eqref{ran} is also the optimal range  for \eqref{car}.  In fact,  by the identity $e^{v\cdot x}(-\Delta)e^{-v\cdot x}= -\Delta +2v\cdot\nabla -|v|^2$ as before,   the estimate \eqref{car} is equivalent to the uniform estimate
\begin{equation}\label{car1}
\|u\|_{L^q(\R^d)} \le C\left\|\left(-\Delta+2v\cdot \nabla - |v|^2 \right)u\right\|_{L^p(\R^d)},
\end{equation}
with a constant $C$ independent of $v\in\R^d$ and $u\in C^\infty_0(\R^d)$.  This estimate is clearly weaker than the uniform Sobolev estimate  \eqref{uni} (see 
Section \ref{proof of theorems}).

Our first result completely characterizes the admissible $p,q$ for which \eqref{car} (equivalently,  \eqref{car1})  holds. 

\begin{thm}\label{carl}  Let $d\ge 3$,  $1<p< q<\infty$. 
The Carleman inequality \eqref{car} holds if and only if 
\begin{equation}\label{range}
\frac 1p-\frac 1q=\frac 2d,   \quad  \frac{d^2-4}{2d(d-1)}  \le \frac1p \le \frac{d+2}{2(d-1)}.
\end{equation}
\end{thm}

The estimate \eqref{car} fails to be uniform for the case $p=1$ or $q=\infty$ since 
the Hardy-Littlewood-Sobolev inequality does not hold either for $p=1$ or $q=\infty$. 
Theorem  \ref{car} gives the estimate \eqref{car} for  $p,q$ which lie outside of the range of admissible $p,q$ for \eqref{uni} (also, see Figure \ref{fig} below). This is in contrast with the non-elliptic differential operators  for which such Carleman estimate is possible  if and only if $p,q$ are admissible exponents for the uniform Sobolev estimate, see Theorem \ref{non-ell-carleman}.  When $d=3, 4$, the second condition in \eqref{range} can be removed. That is to say, for $d=3, 4$, the estimate \eqref{car} holds on the range of $p,q$ where the $L^p$--$L^q$ Hardy-Littlewood-Sobolev inequality is valid.

The estimates \eqref{car} and \eqref{uni} are equivalent to uniform boundedness of  associated multiplier operators. 
Compared with that of the estimate \eqref{uni}, roughly speaking, it may be said  that  the multiplier associated with \eqref{car}  has  singularity on a smaller  set.  This observation is crucial for obtaining the estimate \eqref{car} in an extended range. In order to exploit this  we use the lower $(d-2)$-dimensional restriction-extension estimate associated with the sphere, see Lemma \ref{rextcy} and Section \ref{proof of theorems}.

 By making use of the argument in \cite{KRS}  which shows (weak) unique continuation property for the differential inequality $|\Delta u|\le |Vu|$,  we see that  the extended range of admissible $p,q$ for  the estimate \eqref{car} allows  a larger class of functions for the unique continuation property. 

\begin{cor}\label{unique}  Let $d\ge3$ and $\Omega$ be a connected open set in $\R^d$. Suppose that $u\in W^{2,p}_{loc}(\Omega)$ for some $p>1$ if $d=3,4$,  and $u\in W^{2,\frac{2(d-1)}{d+2}}_{loc}(\Omega)$ if $d\ge5$. Assume that $u$  satisfies $|\Delta u|\le |Vu|$ in $\Omega$ with $V\in L^{d/2}_{loc}(\Omega)$. Then $u$ is identically zero in $\Omega$ whenever $u=0$ in a nonempty open subset of $\Omega$.
\end{cor}

\subsubsection*{Dirac operator in $\mathbb R^2$}  Related to the unique continuation property  for the Dirac operator,  the estimate
\[
 \norm{e^{\lambda \phi(x)}\, \nabla u }_{L^q(\Omega)} \le C\norm{e^{\lambda \phi(x)}\,\Delta u}_{L^p(\Omega)},  \quad u\in C^\infty_0(\Omega)
\]
 which holds uniformly in $\lambda$ 
has been  of interest, where $\Omega$ is an open subset of $\R^d$. It is known that, if $\phi$ is regular, such estimate holds only if $\frac1p-\frac1q\le \frac{2}{3d-2}$ even though $\Omega$ is  bounded (see \cite{BKRS,J, wolff1993}).  A particular example of such estimates is the inequality 
\Be
\label{car-dirac}
 \norm{e^{v\cdot x}\, \nabla u }_{L^q(\R^d)} \le C\norm{e^{v\cdot x}\,\Delta u}_{L^p(\R^d)}. 
 \Ee
By scaling and the Littlewood-Paley decomposition we see that \eqref{car-dirac} is possible only for $p, q$ satisfying the condition 
$\frac1p-\frac1q=\frac1d$.  There is  another necessary condition $\frac1p - \frac1q\ge\frac{2}{d+2}$ (see Lemma \ref{lower} and Theorem \ref{local}). Combining these necessary conditions asserts that the estimate \eqref{car-dirac}  is possible only for $d=2$.  Our method used for the proof of Theorem \ref{carl} also shows that this is indeed the case. 

\begin{thm}\label{dirac} Let $d=2$. Then the estimate \eqref{car-dirac} holds with $C$ independent of $v\in \mathbb R^2$ provided  that $1/p-1/q=1/2$  and $1<p<2$. 
\end{thm}

To our knowledge,  Theorem \ref{dirac} has not been known before and this gives a  complete characterization of  $p,q$ for which \eqref{car-dirac} holds because \eqref{car-dirac}  implies the Hardy-Littlewood-Sobolev inequality. 
Using Theorem \ref{dirac} it is rather straightforward to deduce the Carleman estimates for the Dirac operator $\mathcal D$
\begin{equation} \label{carl_dirac}
\norm{e^{v\cdot x}\, u }_{L^q(\R^2)} \le C\norm{e^{v\cdot x}\, \mathcal D u}_{L^p(\R^2)},\quad u\in C_0^\infty (\R^2;\C^2)
\end{equation}
with $C$ independent of $v\in\R^2$. See Section \ref{dirac_sec} for definition of the  Dirac operator in $\mathbb R^2$. Then, the reflection argument in \cite{KRS} combined with the Kelvin transform for $\mathcal D$ (see Lemma \ref{id-k}) yields the following weak unique continuation property.

\begin{cor}\label{dirac_unique} Let $p>1$ and $\Omega\subset \R^2$ be a connected open set. Assume that $V\in L^2_{loc}(\Omega)$, $u\in W^{1,p}_{loc}(\Omega;\mathbb C^2)$  and 
\begin{equation}\label{diff_dirac}
|\mathcal Du(x)|\le |V(x)u(x)|, \quad x\in\Omega. 
\end{equation}
Then, $u$ is identically zero in $\Omega$ whenever $u$ vanishes in a nonempty open subset of $\Omega$. 
\end{cor}

Again, our contribution here is enlargement of the class of function $u$ for which unique continuation property holds. In higher dimensions $d\ge 3$ unique continuation for \eqref{diff_dirac} with $V\in L^{d}(\R^d)$ can be deduced from  Wolff's result concerning the inequality   
$|\Delta u|\le V|u|+W|\nabla u|$  but it requires a stronger assumption that $u\in W^{2,\frac{2d}{d+2}}_{loc}(\Omega)$, see  \cite[Theorem 1]{wolff1992} for the detail.  For $d=2$, the strong unique continuation property follows from the estimate \cite[Proposition2.6] {wolff1993} for $u\in  W^{2,\,4/3}_{loc}$ (also, see \cite{kim, wolff1990}). 
There is a large body of literature concerning the unique continuation properties for the elliptic and Dirac operators, for example, see  \cite{wolff1993, KoTa} and references therein.

 \subsubsection*{Wave and nonellipitic operators} Let $Q$ be a nondegenerate real quadratic form on $\mathbb{R}^d$, $d\ge 3$, which is given by
\begin{equation}\label{quad}
Q(\xi)=-\xi_1^2-\cdots -\xi_l^2+\xi_{l+1}^2+\cdots+\xi_d^2 ,
\end{equation}
where $1\le l\le d-1$. Let $D=(D_1,\cdots ,D_d)$, $D_j=\frac{1}{i}\frac{\partial}{\partial x_j}$. Appearance of mixed signatures $\pm$ gives rise to a noncompact zero set for  $Q(\xi)$. So, compared with \eqref{uni}, nonelliptic cases exhibit boundedness of different nature. It was shown in  \cite{JKL} (also, see \cite{KRS}) that the uniform Sobolev type estimate associated with $Q$
\Be \label{nonellip-uniform} \| u\|_{q} \le C\|(Q(D)+a\cdot D +b) u\|_p
\Ee
holds with $C$ independent of $(a,b)\in \mathbb C^d\times \mathbb C$,   if and only if  $p,q$ satisfy
\Be\label{nonellip-range} \frac1p-\frac1q=\frac2d\,,\quad    
{p} <\frac{2(d-1)}{d} ,\quad
\frac{2(d-1)}{d-2}<q\,. \footnote{Equivalently, $(1/p,1/q)\in (\fQ(d), \fQ'(d))$ in (I) of Figure \ref{fig}. {For a pair of points $X, Y\in [0,1]\times [0,1]$, we denote by $(X,Y)$ ($[X,Y]$, $[X,Y)$, resp.)  the open (closed, half-open, resp.) line segment connecting $X$ and $Y$. }}
\Ee
Moreover, for $p,q$ satisfying  $(1/p,1/q)=\fQ(d), \fQ'(d)$,  the restricted weak type estimates hold.  The estimate for $p=q'$ is due to Kenig-Ruiz-Sogge \cite{KRS}.  
Main difference from the estimate \eqref{uni} arises in that the zero set of $Q$ is no longer compact and its gaussian curvature vanishes as $|\xi|\to \infty$.  
As before, from  \eqref{nonellip-uniform}  one can deduce the Carleman estimate
\begin{equation}\label{car-non}
\norm{e^{v\cdot x} \, u}_{L^q(\R^d)} \le C\norm{e^{v\cdot x}\, Q(D) u}_{L^p(\R^d)}, 
\end{equation}
which is valid  for $p,q$ satisfying \eqref{nonellip-range} for $d\ge 3$.    
In the following we show that the  optimal range \eqref{nonellip-range} of \eqref{nonellip-uniform} is the same with that of \eqref{car-non}.

\begin{thm} \label{non-ell-carleman} Let $d\ge 3$ and $v\in\mathbb R^d$.  The estimate \eqref{car-non}  holds with $C$ independent of $v\in\mathbb R^d$ if and only if $p,q$ satisfy \eqref{nonellip-range}. 
\end{thm}

\subsubsection*{Heat operator}  We now study the similar  form of Carleman estimate for the heat  operator.  Let $(v,\gamma)\in \mathbb R^{d}\times\mathbb R$ and consider the estimate
\Be \label{heat}  \| e^{(v,\gamma)\cdot(x,t)} F\|_q\le C  \| e^{(v,\gamma)\cdot(x,t)} (\partial_t-\Delta)  F\|_p, \quad F\in C_0^\infty(\R^d\times \R) \Ee
with $C$ independent of $(v,\gamma)\in \mathbb R^{d}\times\mathbb R$. 

\begin{thm}\label{carl-heat}  Let $d\ge 1$ and let $1<p,q<\infty$.  
The estimate \eqref{heat} holds with $C$ independent of $(v,\gamma)$ if  and only if
\begin{equation}\label{hrange}
\frac 1p-\frac 1q=\frac 2{d+2}, \quad \frac1p\ge \frac{d^2+3d-2}{2d(d+2)}, \quad \frac 1q\le\frac{d^2+d+2}{2d(d+2)}.
\end{equation}
\end{thm}

As to be seen later in  Section \ref{section_heat},  unlike the Laplacian  and non-elliptic operators, the estimate \eqref{heat} exhibits a different nature in that the range of  $p,q$ for the uniform bound depends on the direction of  $(v,\gamma)$.  Indeed,    when $0\ge  |v|^2+\gamma$,  \eqref{heat} holds with $C$ independent of $(v,\gamma)$ if and only if $ \frac 1p-\frac 1q=\frac 2{d+2}$.  When $0<  |v|^2+\gamma$, 
\eqref{heat} holds uniformly if and only if \eqref{hrange} is satisfied.

For the time dependent Schr\"odinger operator the Carleman estimate with the weight $e^{(v,\gamma)\cdot(x,t)}$ was obtained by Kenig-Sogge \cite{KS}, and it was extended to a mixed norm setting in Lee-Seo \cite{ls1}, where equivalence between the inhomogeneous Strichartz estimates for the Schr\"odinger equation and the Carleman estimates was established.
As an application the following unique continuation property  can be obtained. 

\begin{cor}\label{heat_unique}
Let $d\ge1$ and assume that $1<p<\frac{d+2}2$ if $d=1,2,$ and $\frac{2d(d+2)}{d^2+5d+2}\le p\le \frac{2d(d+2)}{d^2+3d-2}$ if $d\ge 3.$ Suppose that $V\in L^{\frac{d+2}{2}}(\R^{d+1})$ and $u\in W^{1,p}(\R;W^{2,p}(\R^d))$ satisfies
\begin{equation}\label{diff_ineq}
| (\partial_t-\Delta)u(x,t)|\le |V(x,t)u(x,t)|, \quad (x,t)\in\R^d\times\R.
\end{equation} 
Then, if $u$ vanishes in a half-space in $\R^{d+1}$, $u$ is identically zero on the whole space.
\end{cor}

 The study of the unique continuation property for the heat operator has long history and has also been investigated  by various authors (for example, see \cite{Es,EV,KoTa} and references therein). Especially, the local unique continuation property for a smooth solution $u$ with $V\in L_{loc}^{(d+2)/2}(\R^{d+1})$, $d\ge2,$ was obtained by Sogge \cite{So1}. There have been known various forms of  Carleman estimates for the heat operator,  which is the main tool for the  study of the unique continuation property. But, as far as the authors are aware,  the estimates with the weight $e^{(v,\gamma)\cdot(x,t)}$ such as in Theorem \ref{carl-heat} have not appeared in {the literature} before.

The rest of this paper is organized as follows. In Section 2 we obtain a modification of the restriction-extension estimate \eqref{rext}, which will be used in our argument later. In Section 3 we show main estimates  for related multiplier operators,   
and by making use of them we prove Theorem \ref{carl} and Theorem \ref{dirac}.  
In Section 4 and Section 5, we prove Theorem \ref{non-ell-carleman} and  Theorem \ref{carl-heat}, respectively. In the final section, we provide  proofs of unique continuation results (Corollaries \ref{dirac_unique} and \ref{heat_unique}).

\subsubsection*{Notations} Throughout this paper we write $A \lesssim B $ to denote $A\le CB$ for some constant $C>0$ independent of $A,B>0$. We will write $A \sim B$ to mean $A\lesssim B$ and $B \gtrsim A$.  By $\wh f$\,\, we denote the Fourier transform defined by $\wh f(\xi)=\int e^{-ix\cdot\xi}f(x)dx$, and $\F^{-1}\!f$ is  the inverse Fourier transform.  For a bounded measurable function $m$,  by  $m(D)$ we denote  the Fourier multiplier operator given by  $m(D)f=\F^{-1}(m \wh f \,)$. We  write    
$\xi=(\eta,\tau)\in \mathbb R^{d-1}\times \mathbb R$ and $x=(y,t)\in \mathbb R^{d-1}\times \mathbb R$. By $ \F^{-1}_\tau$, $\F^{-1}_\eta$ we denote the  $1$ and $(d-1)$-dimensional inverse Fourier transforms in $t$ and $y$, respectively.

\section{Preliminaries: A modification of the restriction-extension estimate} 
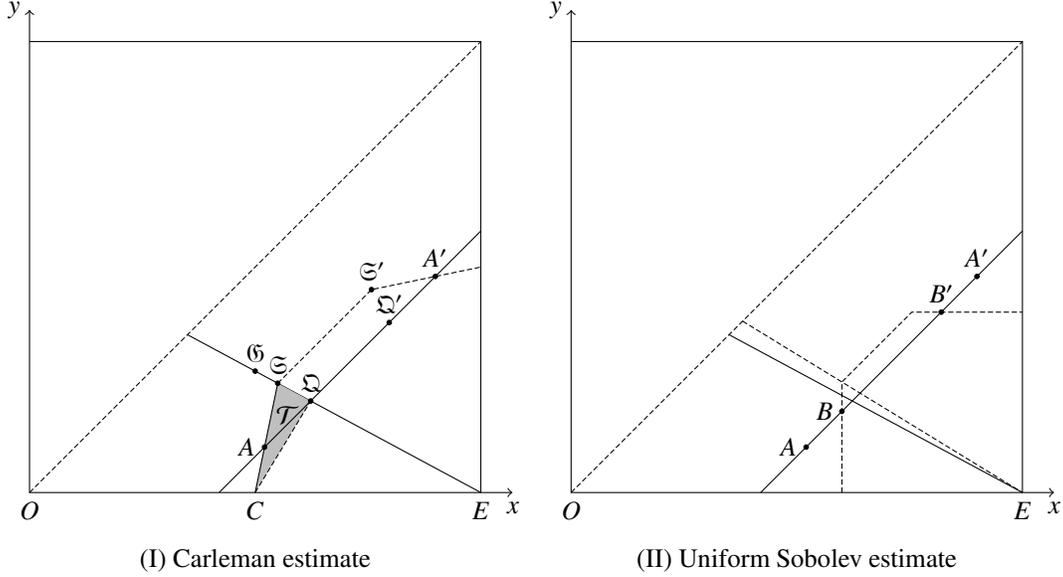
\begin{figure}
\begin{tikzpicture} [scale=0.6]
\begin{scope}
	\draw [<->] (0,10.7)node[left]{$y$}--(0,0) node[below]{$O$}--(10.7,0) node[below]{$x$};
	\draw (0,10) --(10,10)--(10,0) node[below]{$E$};
	\draw (3.5,3.5)--(10,0);
	\path [fill=lightgray] (5,0)--(11/2, 63/26)--(6.23, 2.03); 
	\draw (5.719, 1.7155) node{$\mathcal T$}; 
	\draw (4.2,0)--(10,5.8);
	\draw [dash pattern={on 2pt off 1pt}] (0,0)--(10,10);
	\draw [fill] (5, 35/13) circle [radius=0.05] node[above]{$\mathfrak G$};
	\draw [fill] (6.23, 2.03) circle [radius=0.05] node[above]{$\mathfrak Q$};
	\draw [fill] (10-2.03, 10-6.23) circle [radius=0.05] node[above]{$\mathfrak Q'$};
	\draw [fill] (11/2, 63/26) circle [radius=0.05] node[above]{$\mathfrak S$};
	\draw [fill] (10-63/26, 10-11/2) circle [radius=0.05] node[above]{$\mathfrak S'$};
	\draw (5,0) node[below]{$C$}--(11/2, 63/26); 
	\draw [dash pattern={on 2pt off 1pt}] (5,0)--(6.23,2.03); 
	\draw [dash pattern={on 2pt off 1pt}] (11/2, 63/26)--(10-63/26, 10-11/2)--(10,5);
	\draw [fill] (5.208, 1.008) circle [radius=0.05] node[left]{$A$};
	\draw [fill] (10- 1.008, 10-5.208) circle [radius=0.05] node[above]{$A'$};
	\draw (5,-1.5) node  {(I) Carleman estimate};
\end{scope}
\begin{scope}[shift={(12,0)}]
	\draw [<->] (0,10.7)node[left]{$y$}--(0,0) node[below]{$O$}--(10.7,0) node[below]{$x$};
	\draw (0,10) --(10,10)--(10,0) node[below]{$E$};
	\draw (3.5,3.5)--(10,0);
	\draw [dash pattern={on 2pt off 1pt}] (3.8,3.8)--(10,0);
	\draw [dash pattern={on 2pt off 1pt}] (6,0)--(6, 76/31 )--(10-76/31, 4)--(10, 4);
	\draw (4.2,0)--(10,5.8);
	\draw [dash pattern={on 2pt off 1pt}] (0,0)--(10,10);
	\draw [fill] (5.208, 1.008) circle [radius=0.05] node[left]{$A$};
	\draw [fill] (10- 1.008, 10-5.208) circle [radius=0.05] node[above]{$A'$};
	\draw [fill] (6, 1.8) circle [radius=0.05] node[left]{$B$};
	\draw [fill] (10-1.8, 4) circle [radius=0.05] node[above]{$B'$};
	\draw (5,-1.5) node  {(II) Uniform Sobolev estimate};
\end{scope}
\end{tikzpicture}
\caption{$O=(0,0)$, $A=(\frac{d^2-4}{2d(d-1)},\frac{d-4}{2(d-1)})$, $B=(\frac{d+1}{2d},\frac{d-3}{2d})$, $ C=(\frac12,0)$, $E=(1,0)$, and $\mathfrak{G}=\mathfrak{G}(d),\,\mathfrak{Q}=\mathfrak{Q}(d),\,\mathfrak{S}=\mathfrak{S}(d)$. In the left figure, the closed line segment $AA'$ is  the optimal range of $(1/p,1/q)$ for the Carleman estimate \eqref{car}, which is larger than that for the uniform Sobolev inequality \eqref{uni}, the open line segment $BB'$ in the right.  The line segment  $C\fS$ is contained in the line $  dx-y= {d}/{2}$.}
\label{fig}
\end{figure}

To facilitate statements let us define points  $ \fG(d),$ $\mathfrak Q(d),$  $\fS(d),$ and $\mathfrak Q '(d)$, $ \fS '(d)\in [0,1]\times [0,1]$  by setting 
\begin{align*} 
\fG(d)&=\Bigg(\frac{1}{2},\, \frac{d-2}{2d}\Bigg)\,,   \ \ \ \  \fQ(d)=\Bigg(\frac{d}{2(d-1)}, \,\frac{(d-2)^2}{2(d-1)d}\Bigg)\,,\\ 
&\fS(d)=\Bigg( \frac{d^2+2d-4}{2(d+2)(d-1)}, \, \frac{d(d-2)}{2(d+2)(d-1)}\Bigg )\,,
\end{align*} 
and $(x,y)'=(1-y, 1-x)$. Note that the points $\fG(d)$, $\fQ(d)$, and  $\fS(d)$ are on the line ${y}=\frac{d-2}d(1-x$), see  (I) in   Figure \ref{fig}.

We now recall the following $L^2$ restriction theorem which is known as the Stein-Tomas theorem. 

\begin{thm}[Stein-Tomas  Theorem {\cite{S, T}}]\label{s-t} Let $d\ge 2$. Then, \[  \Normo{\int_{\mathbb S^{d-1}}   e^{ix\cdot\theta} g(\theta) d\theta}_{L^\frac{2(d+1)}{d-1}(\mathbb R^d)} \lesssim \|g\|_{L^2(\mathbb S^{d-1})} .\]  
\end{thm}

The following are the estimates for the restriction-extension operator defined by the sphere which are due to Sogge \cite{So}, Carbery-Soria \cite{CS},  
Bak-McMichael-Oberlin \cite{bmo}, Guti\'{e}rrez \cite{guti} (also, see  \cite{bak} and \cite{CKLS} for more on the  Bochner-Riesz operators of negative order). 

\begin{thm}[Restriction-extension estimates] \label{res-ext}   Let $d\ge 2$. Then, the estimate 
\eqref{rext} 
holds if and only if 
\[    \frac1p-\frac1q\ge \frac{2}{d+1}, \quad    p<\frac{2d}{d+1}, \quad   q >\frac{2d}{d-1}.   \] 
Furthermore, at the critical $p,q$ with $(1/p,1/q)= \fQ(d+1),\, \fQ'(d+1)$,  the restricted weak type  estimate holds. 
\end{thm}

The  estimates \eqref{rext} for $p,q$ satisfying $   \frac1p-\frac1q= \frac{2}{d+1},$  $ p<\frac{2d}{d+1},$   $q >\frac{2d}{d-1}$ were obtained in \cite{bmo} and the restricted weak type endpoint cases were proved in \cite{guti}.  

 Let $L^{p,r}$ denote the Lorentz space.  The following is  a simple modification of Theorem \ref{res-ext}, which we need for the proof of Theorem \ref{carl}.  In what follows we exclusively  use $f$ and $h$ to denote functions on $\mathbb R^d$ and $\mathbb R^{d-1}$, respectively.

\begin{lem}\label{rextcy} Let $d\ge 2$,  $\psi\in C_0^\infty(1/2,2)$, and  $(\frac1p, \frac1q)=\fQ(d)$. Then we have the estimate
\begin{equation}\label{rest-weak}
\Normo{ \int_{1/2}^2 \int_{\mathbb S^{d-2}} \wh f (\phi, \tau) e^{ i ( y \cdot \phi + t \tau ) } d\phi\, \psi (\tau)\, d\tau }_{L^{q,\infty}_{x}(\R^d)} \lesssim \| \psi \|_{C^2} \| f \|_{L^{p,1}(\R^d)}.
\end{equation}
\end{lem}
When  $d=2$  the strong type estimate is trivially true with $p=1$ and $q=\infty$ if we identify $\mathbb S^0=\{-1,1\}$.    Estimate of the strong type is also valid for  $p,q$ satisfying $(1/p,1/q)\in  (\fQ(d), \fQ'(d))$ but these estimates readily follow from  duality,  interpolation  and  easy manipulations. 

{To prove \eqref{rest-weak}, we use the following simple fact which is easy to check: }  For   $1< p<\infty$, 
\begin{equation}
\label{mixed}
 \|f\|_{L_{x}^{p,\infty}}\le \|f\|_{L^p_tL^{p,\infty}_{y}},   \quad \|f\|_{L^p_tL^{p,1}_{y}}\le \|f\|_{L_{x}^{p,1}}, \quad  x=(y,t)\in \mathbb R^{d-1}\times \mathbb R.
 \end{equation}  
 Indeed,   by Fubini's theorem and definition of weak space 
$
 m \{ (y,t): |f(y,t)|>\lambda\} =  \int m_{y}  \{ (y,t): |f(y,t)|>\lambda\}  dt   \le  \lambda^{-p} \int \|f(\cdot,t)\|_{L^{p,\infty}_{y}}^p  dt  =\lambda^{-p}  \|f\|_{L^p_tL^{p,\infty}_{y}}^p.
 $ This gives the first inequality and the latter inequality in \eqref{mixed} follows from duality.

\begin{proof}  Let   $(\frac1p, \frac1q)=\fQ(d)$, and from Theorem \ref{res-ext}  recall the following estimate for the restriction-extension operator defined by the $(d-2)$-sphere:
\begin{equation}\label{rext1}
\Normo{\int_{\mathbb S^{d-2}} \wh h(\phi)e^{iy\cdot\phi}d\phi}_{L^{q,\infty}_y(\R^{d-1})} \lesssim \|h\|_{L^{p,1}(\R^{d-1})}. 
\end{equation}
Since
\[
 \iint_{\mathbb S^{d-2}} \wh f (\phi, \tau) e^{ i ( y \cdot \phi + t \tau ) } d\phi\, \psi (\tau) d\tau = 2\pi  \int \F^{-1}_\tau(\psi )(t-s) \int_{\mathbb S^{d-2}} \F ( f (\cdot, s) ) (\phi) e^{iy\cdot\phi} d\phi \,ds,
\]
 the first estimate in \eqref{mixed} and Minkowski's inequality give 
\begin{align*}
\Normo{ \iint_{\mathbb S^{d-2}} \wh f (\phi, \tau) e^{ i ( y \cdot \phi + t \tau ) } d\phi\, \psi (\tau) d\tau }_{L_{y,t}^{q,\infty}} 
& \lesssim  \Normo{ \int  \F^{-1}_\tau(\psi )(t-s) \int_{\mathbb S^{d-2}} \F ( f (\cdot, s) ) (\phi) e^{iy\cdot\phi} d\phi \,ds}_{L^{q}_t  L_{y}^{q,\infty}} 
\\
&\le    \Normo{ \int \left| \F^{-1}_\tau(\psi )(t-s) \right| \normo{ \int_{\mathbb S^{d-2}} \F ( f (\cdot, s) ) (\phi) e^{iy\cdot\phi} d\phi }_{L_{y}^{q,\infty}} ds}_{L^{q}_t}.
\end{align*}
Using the restriction-extension estimate \eqref{rext1} and  Young's convolution inequality, we get
\begin{align*}
\Normo{ \iint_{\mathbb S^{d-2}} \wh f (\phi, \tau) e^{ i ( y \cdot \phi + t \tau ) } d\phi \psi (\tau) d\tau }_{L_{x}^{q,\infty}} 
&\lesssim   \Normo{ \int \left|\F^{-1}(\psi )(t-s) \right| \Normo{  f (\cdot, s)}_{L_{y}^{p,1}(\R^{d-1})} ds}_{L^{q}_t} 
\\
&\le \Normo{\F^{-1}(\psi )}_r   \norm{ f}_{L_t^p  L_{y}^{p,1}}  
\lesssim \Normo{\F^{-1}(\psi )}_r \norm{f}_{L^{p,1}},
\end{align*}
where $1/r=1/q-1/p+1\in [0,1]$. Since $\supp \psi \subset [1/2,2]$ it is clear that $\|\wh\psi\|_\infty \le \|\psi\|_1 \lesssim \|\psi\|_\infty$, and  $\|\wh\psi\|_1\lesssim \|\psi\|_\infty+\|\psi''\|_\infty$ from integration by parts. Therefore we have $\|\F^{-1}(\psi )\|_r\lesssim \|\psi\|_{C^2}$ for any $r\in[1,\infty]$.
\end{proof}

\section{Estimates for $\Delta$: Proofs of Theorem \ref{carl} and Theorem \ref{dirac}} \label{proof of theorems}

As mentioned before,  for  the estimate \eqref{car} it is enough to show \eqref{car1}. By taking the Fourier transform we observe that  \eqref{car1} is equivalent to 
\[
\normo{\F^{-1}\left(\frac{\wh f(\xi)}{|\xi|^2+2iv\cdot\xi -|v|^2} \right)}_{L^q(\R^d)} \le C\|f\|_{L^p(\R^d)},
\]
with the constant $C$ independent of $v\in\R^d$.  Thanks to the homogeneity condition $1/p-1/q=2/d$ and rescaling $(\xi\to |v|\xi)$,  we may assume that $|v|=1$. Also, by rotation, we may further  assume that $v=e_d=(0,\cdots,0,1)$. Thus, for  the  inequality \eqref{car1} it is sufficient to show the following {estimate}:
\begin{equation}\label{m}
\normo{\F^{-1}\left(\frac{\wh f(\xi)}{|\xi|^2+2i\xi_d -1} \right)}_{L^q(\R^d)} \le C\|f\|_{L^p(\R^d)}.
\end{equation}

 Let $\chi$ be a smooth function such that $\chi=1$ on $B(0,3/2)$ and supported in  $B(0,2)$. If $|\xi|\ge 3/2$, $A(\xi)=\frac {|\xi|^2}{|\xi|^2+2i\xi_d -1}$ satisfies $\partial_\xi^\alpha A(\xi)=O(|\xi|^{-|\alpha|})$. Thus by Mihlin's multiplier theorem and the Hardy-Littlewood-Sobolev inequality we see 
 \[ \normo{ \F^{-1} \left(\frac{\wh f(\xi)(1-\chi(\xi))}{|\xi|^2+2i\xi_d -1} \right)}_q\lesssim \normo{\F^{-1}\left(\frac{\wh f(\xi)}{|\xi|^2} \right)}_q\lesssim \|f\|_p\] for $1<p,q<\infty$ satisfying $1/p-1/q=2/d$. Thus, for the proof of Theorem \ref{carl} it is sufficient to consider the multiplier operator given by \[ m_{\Delta}(\xi)=\frac 1{|\xi|^2+2i\xi_d -1}  \chi(\xi).\]
For technical convenience we assume that $\chi$ is radial in the first $(d-1)$-variables and write {$\chi(\rho\theta, \tau)=\chi_\circ(\rho,\tau)$ when $\theta\in\mathbb S^{d-2}$, $\rho>0$.}

{
\begin{thm}\label{local} Let $1<p, q<\infty$, and  
 $d\ge 2$. We have the estimate
\Be 
\label{m}
\normo{\F^{-1}\left(m_{\Delta}(\xi)\wh f(\xi) \right)}_{L^q(\R^d)} \le C\|f\|_{L^p(\R^d)}
\Ee 
if and only if $p,$ $q$ satisfy 
\begin{align}
 \label{nec-cond}
 \frac1p-\frac1q\ge \frac{2}{d+2}, \quad
  \frac dp-\frac1q\ge \frac{d}{2}, \quad  \frac{d-2}2 \ge  \frac dq-\frac1p\,.
 \end{align}  
\end{thm}
}
The pairs of $(p,q)$ satisfying \eqref{nec-cond}  are those of $(p,q)$ for which $(1/p,1/q)$ is contained in  the closed pentagon with vertices $C, \fS, \fS', C', E$ in (I) of Figure \ref{fig}.

The key observation is that the multiplier $m_\Delta$  essentially localizes $O(\tau)$-neighborhood of $\mathbb S^{d-2}$ in $\mathbb R^{d-1}$ if $|\xi_d|\sim \tau$. This leads us to 
decompose dyadically the multiplier $m_\Delta$ along $\xi_d$ and use the restriction-extension estimates associated with  the $(d-2)$-dimensional submanifold $ \{\xi\in\R^d: |\xi|=1, \, \xi_d=0\}.$ However, it is important to obtain estimates which are sharp in terms of the value of $\tau$.

For any bounded measurable function  $G$ we define the multiplier operator norm from $L^{p,r}$ to $L^{q,s}$ by 
\[   \mathfrak M_{(p,r), (q,s)} [G]  =    \sup \left\{ \normb{  \finv{ G\widehat f\,}   }_{L^{q,s}(\R^d)} : \|f\|_{L^{p,r}(\R^d)} =1, \  f\in\mathcal S(\mathbb R^d) \right\}.  \] 
When $p=r$ and $q=s$ we use a simplified notation $\mathfrak M_{p, q}$ for 
$\mathfrak M_{(p,p), (q,q)}$.

We fix a smooth cut-off function $\psi\in C_0^\infty( 1/2,2 )$ such that $\sum_{j=-\infty}^\infty \psi(t /2^j)=1$ for $t> 0$ and let $\mathbb D$ be the set of positive dyadic numbers contained in $[0,2]$.  We decompose the multiplier 
\Be\label{1st-decomp}
m_\Delta(\xi) =  \sum_{\varepsilon\in \dn} m_\varepsilon(\xi ) :=   \sum_{\varepsilon\in \dn}      \psi(|\tau| /\varepsilon)  m_\Delta(\xi) , \quad \xi=(\eta,\tau)\in\R^{d-1}\times\R.
\Ee
If $\varepsilon$ is large, that is to say $\varepsilon \ge \varepsilon_0$ for some  $\varepsilon_0>0$, then $\sum_{\varepsilon\ge \varepsilon_0 }m_\varepsilon$ is a smooth multiplier with compact support and $\| \sum_{\varepsilon\ge \varepsilon_0 }m_\varepsilon \|_{C^N}\le C(\varepsilon_0, N)$ for any $N$. Hence, 
\[ \Mpqq{p}{q}{\sum_{\varepsilon\ge \varepsilon_0}m_\varepsilon}\le C\]  for any $1\le p\le q\le \infty$. In what follows we may assume that $\varepsilon \le \varepsilon_0$ for $\varepsilon_0>0$ small enough. By the similar argument, it is easy to see that $\mathfrak M _{p,q} [ \varphi m_\Delta ]\le C$ for any $1\le p\le q\le \infty$, if $\varphi$ is a smooth function which vanishes on the set $\{\xi: ||\eta|^2-1|\le \varepsilon_0\}$. Hence, from now on, we may also assume that $m_\Delta =(1-\varphi)m_\Delta$.

\begin{lem}\label{delta-local} Let $d\ge 2$.   Suppose  $p, q$ satisfy that  $(\frac1p,\frac1q)$ is contained in the triangle $\mathcal T$ with vertices $\mathfrak S(d), \mathfrak Q(d), (\frac12,0)$ from which the line segment $[\mathfrak Q(d), (\frac12,0))$ is removed. (See (I) in Figure \ref{fig}). Then we have 
\begin{equation}\label{tri}   \mathfrak M_{p, q} [m_\varepsilon]\lesssim \varepsilon^{\frac dp-\frac1q-\frac d2}.
\end{equation}
In particular, if $(\frac1p,\frac1q)\in [\mathfrak S(d), \mathfrak Q(d)\big)$, 
we have the estimate
\begin{equation}\label{line-dg}   \mathfrak M_{p, q} [m_\varepsilon]\lesssim \varepsilon^{-1+\frac{d+2}2(\frac1p-\frac1q)}.
\end{equation}
\end{lem}

\newcommand{\eps}{\varepsilon} 
The range of $p,q$ can be extended by duality and interpolation with other easy estimates. But the resulting estimates are irrelevant to our purpose.  Before we provide  proof of Lemma \ref{delta-local}  we show that the range \eqref{nec-cond} in Theorem \ref{local} is optimal.  Meanwhile, we also see the bounds in Lemma \ref{delta-local} are sharp.

\begin{lem}\label{lower}Let $d\ge2$, $0<\varepsilon\ll1$, and $1\le p,q\le \infty$. Then we have 
\begin{align} 
\label{lower1}
 \mathfrak M_{p, q} [m_\varepsilon]&\gtrsim \varepsilon^{-1+\frac{d+2}2(\frac1p-\frac1q)},  \\
\label{lower2}
\mathfrak M_{p, q} [m_\varepsilon]&\gtrsim \varepsilon^{\frac dp-\frac1q-\frac d2}.
\end{align}
\end{lem}

\begin{proof} 
To begin with, we observe that 
$\mpq{m_\varepsilon}=\mpq{\overline m_\varepsilon}$. Thus, by the triangle inequality  
\Be \label{real}\mpq{m_\varepsilon}\ge  \mpq{\im (m_\varepsilon)}.  \Ee

Via scalng $\tau\to \eps\tau$ we have that,   for any $\varepsilon>0$, 
\begin{equation}\label{nec}
\mpq{\im (m_\varepsilon)}= \varepsilon^{\frac 1p-\frac 1q} \mpq{\rho_\varepsilon} = \varepsilon^{\frac 1p-\frac 1q} \sup_{0\neq f \in L^p(\R^d)} \frac{ \|\rho_\varepsilon(D)f\|_q }{\|f\|_p},
\end{equation}
where  
\[
\rho_\varepsilon (\eta, \tau) :=  -\im (m_\varepsilon) (\eta, \varepsilon\tau) = \frac{2\varepsilon\tau \chi(\eta, \varepsilon\tau)}{(|\eta|^2-1+\varepsilon^2\tau^2)^2+4\varepsilon^2\tau^2}\psi(|\tau|), \quad (\eta, \tau) \in \R^{d-1}\times \R\,.
\]

Let us set
\[
\widetilde{\rho_\varepsilon }(D)f(x)=\iint \frac{ 2\varepsilon\tau \chi(\eta-{e_{d-1}}, \varepsilon\tau)}{(|\eta|^2-2\eta_{d-1}+\varepsilon^2\ta^2)^2+4\varepsilon^2\ta^2} \psi(|\ta|)\wh f(\eta,\ta)e^{i(y\cdot\eta+t\ta)}d\eta\, d\ta.
\]
The above multiplier only differs from the previous one  by  translation $\eta\to \eta -e_{d-1}$. So, it is clear that ${ \|\rho_\varepsilon(D)f\|_q }/{\|f\|_p}={ \|\widetilde {\rho_\varepsilon}(D)f\|_q }/{\|f\|_p}$.  We pick a function $\phi\in C_0^\infty(\R)$ such that $\phi\ge0$, $\phi(0)=1$ and $\supp\phi \subset [-1,1]$, and  put
\[
\wh f(\eta, \tau)=\psi(\tau) \phi \left(\frac{\eta_{d-1}}\varepsilon \right)\prod_{j=1}^{d-2}\phi\left(\frac{10d\,\eta_j}{\sqrt{\varepsilon}}\right), \quad  \eta=(\eta_1,\dots, \eta_{d-1} ).
\]
Set $A=\{x\in\R^d: |x_j|\le c \varepsilon^{-1/2}, 1\le j\le d-2, \, |x_{d-1}|\le c \varepsilon^{-1}, \, |x_d|\le c\}
$ 
for some $0< c\le 1$. 
Then, it is easy to see that  if $c>0$ is small enough and $x\in A$
\[
\left|\widetilde{\rho_\varepsilon}(D)f(x)\right|\gtrsim \varepsilon^{-1}\varepsilon^{\frac {d-2}2+1}.
\]
Consequently, $\|\widetilde{\rho_\varepsilon}(D)f\|_q\gtrsim \varepsilon^{\frac {d-2}2} (\varepsilon^{-\frac {d}2})^\frac 1q $. Thus, from \eqref{nec} we have
\[
\mpq{\im (m_\varepsilon)} 
\ge 
\varepsilon^{\frac 1p-\frac 1q} 
\frac{ \|\widetilde{\rho_\varepsilon}(D)f\|_q }{\|f\|_p} 
\gtrsim 
\frac{\varepsilon^{\frac 1p-\frac 1q}\varepsilon^{\frac {d-2}2} (\varepsilon^{-\frac {d}2})^\frac 1q }{\varepsilon^{\frac d2}(\varepsilon^{-\frac {d}2})^\frac 1p } =\varepsilon^{-1+\frac{d+2}2(\frac1p-\frac1q)}.
\]
Thus, by \eqref{real}  we get \eqref{lower1}. 

We now turn to  \eqref{lower2}. By \eqref{real} and \eqref{nec} it is enough to show that  $\mathfrak M_{p, q} [\rho_\varepsilon]\gtrsim \varepsilon^{\frac {d-1}p-\frac d2}$ which, by duality, is equivalent to the following: 
\Be\label{lower-} \mpq{{\rho_\varepsilon }}\gtrsim  \varepsilon^{\frac{d-2}2-\frac{d-1}q}.\Ee

We choose a special $f$ such that \[\wh f (\eta,\tau) =\chi_0(|\eta|)\chi_0(\tau),\] 
where $\chi_0$ is a smooth function such that  $0\le \chi_0\le 1$, $\chi_0(\rho)=1$ if $|\rho-1|\le \delta_\circ$, and $\chi_0(\rho)=0$ if $|\rho-1|\ge 2\delta_\circ$.  Here $\delta_\circ>0$ is a fixed small number.  
Clearly, we have  
\[\mpq{{\rho_\varepsilon }}\gtrsim \|{\mathcal F^{-1}(\rho_\varepsilon \wh f\,)} \|_q. \] 
For $0<\varepsilon \ll1$, $\chi(\eta, \varepsilon\tau)\equiv 1$ on $\supp \wh f$\,, hence
\begin{align*}
\mathcal F^{-1} (\rho_\varepsilon \wh f \,)(y, t)=(2\pi)^{-d}\iint \frac{2\varepsilon \tau \,\chi_0(\rho) \chi_0(\tau)}{(\rho^2-1+\varepsilon^2\tau^2)^2+4\varepsilon^2\tau^2}\psi(\tau) e^{i\tau t} \left(\int_{\mathbb S^{d-2}}  e^{i\rho \phi\cdot y} d\phi\right)\rho^{d-2}d\tau d\rho.
\end{align*}

Since $\int_{\mathbb S^{d-2}} e^{iy\cdot\phi} d\phi= C_d|y|^{-\frac{d-2}2} \cos(|y|-\frac{d-2}{4}\pi)+O(|y|^{-\frac{d}2})$ and $\cos(\alpha+\beta)=\cos\alpha\cos\beta-\sin\alpha\sin\beta$, we see that  
\[ \int_{\mathbb S^{d-2}}  e^{i\rho \phi\cdot y} d\phi=  C_d|\rho y|^{-\frac{d-2}2} \cos\left(|y|-\frac{d-2}{4}\pi\right) \cos(|y|(1-\rho))+  O\B(\delta_\circ  |\rho y|^{-\frac{d-2}2}\B)\]
{whenever $|y|-\frac{d-2}4\pi \in 2\pi \Z +[ -\delta_\circ,\delta_\circ]$ and $|y|\gtrsim\delta_\circ^{-1}$.   }
Putting this back into the above equation, we have that, for a sufficiently small $\delta_\circ>0$, 
%
\[
| \mathcal F^{-1} (\rho_\varepsilon \wh f \,) (y, t ) | \gtrsim   \varepsilon^{\frac{d-2}2}\cos\left(|y|-\frac{d-2}{4}\pi\right)
\]
provided that  $|y|-\frac{d-2}4\pi\in [c\varepsilon^{-1}, 2c\varepsilon^{-1}]\,\cap\,[2m\pi -\delta_\circ,2m\pi +\delta_\circ]$ for any $m\in \Z$  and $|t|\le c$ with a sufficiently small $c>0$. 
Hence, it follows that 
\[
\|\mathcal F^{-1}(\rho_\varepsilon \wh f\,) \|_q\gtrsim \varepsilon^{\frac{d-2}2-\frac{d-1}q} .\] 
This gives the desired lower bound \eqref{lower-}. 
\end{proof}


We now prove Lemma \ref{delta-local}. 
\begin{proof}[Proof of Lemma \ref{delta-local}
]
We set 
\[
\wt m_\varepsilon (\eta, \tau) :=   m_\varepsilon (\eta, \varepsilon \tau).
\]
By scaling $\tau\to\varepsilon\tau$ it is easy to check that
$
\mpq{m_\varepsilon} =  \varepsilon^{\frac1p-\frac 1q} \mpq{\wt m_\varepsilon} .
$
Thus it is sufficient to show that 
\[
\mpq{\wt m_\varepsilon} \lesssim  \varepsilon^{\frac {d-1}p-\frac d2}
\] whenever $(\frac1p, \frac1q)\in \mathcal T$.

Let us set
\[
\psi_0(\tau)=\sum_{j\le0}   \psi(2^{-j}|\tau|), \quad  \psi_k(\tau)=\psi(2^{-k}|\tau|), \quad  k\ge 1.
\]
 Using this we decompose the multiplier $\wt m_{\varepsilon}$ dyadically in $\eta$ off the $(d-2)$-sphere $|\eta|=1$. For $k\ge 0$, set
\[
\wt m_{\varepsilon, k}(\eta, \tau):=\psi_k \left ( \frac{|\eta|^2-1}{\varepsilon} \right) \wt m_{\varepsilon}(\eta, \tau)=   
\psi_k \left ( \frac{|\eta|^2-1}{\varepsilon} \right)\frac { \psi(|\tau|) \chi(\eta,\epsilon \tau)} {|\eta|^2+\varepsilon^2\tau^2+2i\varepsilon \tau -1}\,  
\]
so that 
\begin{equation}\label{k-sum}  \wt m_\varepsilon= \sum_{k\ge 0}  \wt m_{\varepsilon, k}\, .\end{equation}

The estimate for $\mpqq 2\infty{\wt m_{\varepsilon,k}}$ easily follows from Plancherel's theorem and the Cauchy-Schwarz inequality since $\|\wt m_{\varepsilon,k} \|_2 \lesssim (2^k\varepsilon)^{-1/2}$. This yields $\mpqq 2\infty{\wt m_{\varepsilon}}\lesssim \varepsilon^{-1/2}$. Thus  by interpolation it is sufficient  to show,  for $(\frac1p,\frac1q)\in [\mathfrak S(d), \mathfrak Q(d)\big)$, 
\begin{equation}\label{scaling}   \mathfrak M_{p, q} [\wt m_\varepsilon]\lesssim \varepsilon^{-1+\frac{d}2(\frac1p-\frac1q)}.
\end{equation}  
We first show this when $d\ge 3$.   The case $d=2$ is much easier. 

Firstly  we  show that
\begin{equation}\label{decomp-k}
\mpqq 2{\frac{2d}{d-2}}{\wt m_{\varepsilon, k}}  \lesssim \varepsilon^{-1/2} 2^{-k/2}.
\end{equation}
Using the spherical coordinates we notice that 
\begin{align*}
{\wt m_{\varepsilon, k}} (D)f(x)= \frac 1{(2\pi)^d}\int_0^\infty \psi_k \left ( \frac{r^2-1}{\varepsilon} \right )
         \iint_{\mathbb S^{d-2}} \wt m_\varepsilon(r\phi, \tau)\wh f(r\phi , \tau) e^{i(r\phi\cdot y +\tau t)} d\phi\, d\tau\,  r^{d-2}dr.
        \end{align*}
Let us fix a smooth cut-off function $\wt \psi\in C_0^\infty (-4,4)$ such that $\wt\psi=1$ on $[-2,2]$. For every $k> 0$, by Minkowski's inequality it follows that 
\[
\| \wt m_{\varepsilon, k} (D)f \|_{\frac{2d}{d-2}}\lesssim \int \left |\wt \psi \left( \frac{r^2-1}{2^k\varepsilon}\right) \right| \left|\psi_k
   \left ( \frac{ r^2-1}{\varepsilon} \right ) \right| \Normo{\iint_{\mathbb S^{d-2}} \wt m_\varepsilon  (r\phi, \tau) \wh f(r\phi, \tau)e^{i(r\phi\cdot y +\tau t)} d\phi d\tau}_{L^\frac{2d}{d-2}_{y, t}(\R^d)} dr.
 \]
By Hausdorff-Young's and Minkowski's inequalities we have that, for $r
\sim 1$, 
\[ 
\Normo{\iint_{\mathbb S^{d-2}} \wt m_\varepsilon  (r\phi, \tau) \wh f(r\phi, \tau)e^{i(r\phi\cdot y +\tau t)} d\phi d\tau}_{L^\frac{2d}{d-2}_{y, t}(\R^d)}
\lesssim 
\bigg\| \left \| \int_{\mathbb S^{d-2}} \wt m_\varepsilon  (r\phi, \tau) \wh f(r\phi, \tau)e^{ir\phi\cdot y} d\phi \right\|_{L^\frac{2d}{d-2}_{y}} \bigg\|_{L^\frac{2d}{d+2}_\tau} {.}
\]
Using  the $L^2$-Fourier extension (adjoint  restriction)  estimate  from the sphere $\mathbb S^{d-2}$ (Theorem \ref{s-t}) and H\"older's inequality (in $\tau$) we have that 
\[ \Normo{\iint_{\mathbb S^{d-2}} \wt m_\varepsilon  (r\phi, \tau) \wh f(r\phi, \tau)e^{i(r\phi\cdot y +\tau t)} d\phi d\tau}_{L^\frac{2d}{d-2}_{y, t}(\R^d)}
\lesssim 
\Big\| \wt m_\varepsilon  (r\phi, \tau) \wh f(r\phi, \tau) \Big\|_{{L^2_{\phi,\tau}(\mathbb{S}^{d-2}\times\R)}}  {.}\]
Putting this in the above inequality and using  the Cauchy-Schwarz inequality we see that
\begin{align*}
\| \wt m_{\varepsilon, k} (D)f \|_{\frac{2d}{d-2}}\lesssim  C(2^k\varepsilon)^{1/2} \left( \int \left| \psi_k \left ( \frac{r^2-1}{\varepsilon} \right ) \right|^2\normo{\wt m_\varepsilon(r\phi, \tau) \wh f(r\phi, \tau)}_{L^2_{\phi,\tau}(\mathbb{S}^{d-2}\times\R)}^2 dr \right)^{1/2}
\end{align*}
because $|r-1|\le\varepsilon_0$. Since $|r^2-1|\approx 2^k\varepsilon$ and $\varepsilon\ll 1$, $|\wt m_\varepsilon(r\phi, \tau)| \lesssim (2^k\varepsilon)^{-1}$. Thus, by Plancherel's theorem we have 
\[\normo{  {\wt m_{\varepsilon, k}} (D)f }_{\frac{2d}{d-2}}\lesssim 
(2^k\varepsilon)^{-1/2}  \normo{\psi_k \left ( \frac{|\eta|^2-1}{\varepsilon} \right ) \wh f(\eta, \tau)}_{L^2_{\eta, \tau}(\R^d)}
\lesssim   (2^k\varepsilon)^{-1/2}  \|f\|_2.\]
So, we get \eqref{decomp-k} when $k>0$. The same argument also {works for the case $k=0$.}

Next, we show  
\begin{equation}\label{decomp-kk}
\mpqq {(p,1)}{(q,\infty)}{\wt m_{\varepsilon, k}}  \lesssim  1,  \quad  (1/p, 1/q)=\fQ(d).
\end{equation}

As before note that 
\[
\wt m_{\varepsilon, k}(D)f (y, t) = \frac 1{(2\pi)^d}\int_0^\infty \wt \psi \left ( \frac{\rho^2-1}{2^k\varepsilon} \right ) \int  \chi_{\varepsilon,k}(\rho, \tau) \int_{\mathbb S^{d-2}} \wh f(\rho\phi , \tau) e^{i(\rho\phi\cdot y +\tau t)} d\phi d\tau \rho^{d-2}d\rho
\]
for $k=0, 1,2,\cdots ,$ where
\[
\chi_{\varepsilon, k}(\rho, \tau)= \psi_k \left ( \frac{\rho^2-1}{\varepsilon} \right)\frac {\psi(|\tau|) \chi_\circ(\rho, \epsilon \tau)} {\rho^2+\varepsilon^2\tau^2+2i\varepsilon \tau -1} .
\]
By Minkowski's and H\"older's inequalities, we see that 
\[\normo{\wt m_{\varepsilon, k}(D)f }_{q,\infty} 
\lesssim (2^k\varepsilon) \sup_{\rho: |\rho^2-1|\le \varepsilon_0} \Normo { \int \int_{\mathbb S^{d-2}} \chi_{\varepsilon,k}(\rho, \tau) \wh f(\rho\phi , \tau) 
e^{i(\rho\phi\cdot y +\tau t)} d\phi\,d\tau }_{L^{q,\infty}(\R^d)}. \] 
Since   $\sup_{\rho: |\rho^2-1|\le \varepsilon_0} \|\chi_{\varepsilon,k}(\rho, \cdot)\|_{C^2}\lesssim 1/(2^k\varepsilon)$, 
by making use of Lemma \ref{rextcy} and taking the support of the multiplier $\chi_{\eps,k}$
into account we have
\begin{align*}
\normo{\wt m_{\varepsilon, k}(D)f }_{q,\infty} 
&\lesssim (2^k\varepsilon) \left( \sup_{\rho: |\rho^2-1|\le \varepsilon_0} \normo{\chi_{\varepsilon,k}(\rho, \cdot)}_{C^2} \right) \|f\|_{L^{p,1}(\R^d)} \lesssim \|f\|_{p,1}.
\end{align*} 
This gives \eqref{decomp-kk}.

Now, interpolation between  \eqref{decomp-k} and \eqref{decomp-kk} yields 
\[  \mpq {\wt m_{\varepsilon,k} }\lesssim \varepsilon^{-1+\frac{d}2(\frac1p-\frac1q)}  2^{-\alpha k} \]
for some $\alpha>0$ whenever $(\frac1p, \frac1q)\in [\fG(d), \fQ(d))$. Combining this with \eqref{k-sum} {and} summation along $k$ {give} \eqref{scaling} for $p,q$ satisfying  $(\frac1p, \frac1q)\in [\fG(d), \fQ(d))$. This completes 
the proof of \eqref{scaling} when $d\ge 3$.

For $d=2$ the same argument also works if we identify $\mathbb S^{0}=\{-1,1\}$.  The estimates \eqref{decomp-k} and \eqref{decomp-kk} 
correspond to $ \mpqq {2}{\infty}{\wt m_{\varepsilon, k}}  \lesssim \varepsilon^{-1/2} 2^{-k/2}$ and  $\mpqq {1}{\infty}{\wt m_{\varepsilon, k}}  \lesssim  1$, respectively. The first follows from the Cauchy-Schwarz inequality and Plancherel's theorem, and the latter is clear since 
$\|\widehat f\|_\infty\le \|f\|_1$.  Interpolation and summation along $k$ give the desired estimates. 
\end{proof}

 Once we have Lemma \ref{delta-local}{,}
 the proof of Theorem \ref{local} is rather routine.

 \begin{proof} [Proof of Theorem \ref{local}] We first prove the necessity part.  
 We only need to show that \eqref{m} implies the first and the second inequalities in \eqref{nec-cond}. The third one follows from the second via duality.     
   From the assumption $\mpq{m_\Delta}<\infty$ it follows that $\mpq{m_\varepsilon}<C$ independently of $\varepsilon\ll 1$.  By Lemma \ref{lower}
   \[ C\gtrsim \varepsilon^{-1+\frac{d+2}2(\frac1p-\frac1q)},   \quad
 C\gtrsim \varepsilon^{\frac dp-\frac1q-\frac d2}\]
for   $\varepsilon \ll 1$. Hence, considering the limiting case $\varepsilon\to 0$ yields the first and the second inequalities in \eqref{nec-cond}. 
 
 We now turn to the sufficiency part of Theorem \ref{local}.  Since $m_\Delta$ has compact support, 
 by duality and  interpolation it is sufficient to show that \eqref{m} holds for 
 $(1/p,1/q)\in [\fS(d), (\frac12,0))$.  Indeed, duality gives the estimates for  $(1/p,1/q)\in [\fS(d)', (1,\frac12))$ and interpolation between these estimates  gives the desired boundedness for  $(1/p,1/q)$ on the trapezoid with vertices $\fS(d)',\fS(d), (\frac12,0), (1,\frac12)$ from which the line segment  $[(\frac12,0), (1,\frac12)]$ is removed. Since $m_\Delta$ is compactly supported, we have \eqref{m} for all $(1/p,1/q)$ satisfying \eqref{nec-cond}.  
 
 Let $\overline{\psi}\in C_0^\infty((-4,-1/4)\cup(1/4,4))$ such that $\ovpsi(\tau)\psi(|\tau|)=\psi(|\tau|)$. Then from the Littlewood-Paley theory (\cite{St-singular, S}) we have, for $1<r<\infty$,  
  \Be \label{L-P}
 \|  f\|_{L^r(\mathbb R^d)} \sim   \B\| \B(\sum_{ j\in \Z}\B|\ovpsi\Big(\frac{D_d}{2^{j}}\Big) f\B|^2\B)^\frac12\B\|_{L^r(\mathbb R^d).}\footnote{It can be easily seen by making use of $ \sum_{j\in \Z}\ovpsi^2(2^{-j}t)\sim 1$ and the standard argument for the Littlewood-Paley inequality. In fact, by Fubini's theorem it is enough to show this for functions on $\mathbb R$. }
 \Ee
Since $\ovpsi\Big(\frac{D_d}{\eps}\Big) m_\eps(D)=m_\eps(D)$ and $q\ge 2$, by \eqref{1st-decomp}, the Littlewood-Paley inequality and  Minkowski's inequality, we have 
\begin{align}
\|m_\Delta(D) f\|_q\lesssim  \B\| \B(\sum_{\eps\in \dn}\B|\ovpsi\Big(\frac{D_d}{\eps}\Big) m_\eps(D) f \B|^2\B)^\frac12\B\|_{q}
 \lesssim 
\B(  \sum_{\eps\in \dn}\B\| m_\eps(D)  \ovpsi\Big(\frac{D_d}{\eps}\Big)  f\B\|^2_q \B)^\frac12.
\end{align}
From this and Lemma \ref{delta-local},  we see that if 
$(\frac1p,\frac1q)\in \mathcal T$, 
\begin{align*}
\|m_\Delta(D) f\|_q\lesssim 
\B( \sum_{\eps\in \dn}   \varepsilon^{\frac dp-\frac1q-\frac d2} \B\| \ovpsi\Big(\frac{D_d}{\eps}\Big)   f\B\|^2_p \B)^\frac12 .
\end{align*}
Note that $[\fS(d), (\frac12,0))\subset \mathcal T$, and  $\frac dp-\frac1q-\frac d2=0$ if $(\frac 1p,\frac1q)\in [\fS(d), (\frac12,0))$. Thus,  by Minkowski's inequality with $p\le 2$ we have 
\[\|m_\Delta (D) f\|_q 
\lesssim  \B\|   \B( \sum_{\eps\in \dn}  \B|\ovpsi\Big(\frac{D_d}{\eps}\Big)   f \B|^2\B)^\frac12\B\|_p\lesssim \|f\|_p\]
whenever  $(\frac 1p,\frac1q)\in [\fS(d), (\frac12,0))$.

Finally we consider the case $d=2$. We note $\fQ(2)=(1,0)$ and $\fS(2)=(\frac12,0)$. Thus we have \eqref{line-dg} for $p,q$ satisfying $1/p-1/q\ge 1/2$.  We may repeat the above argument and we obtain \eqref{m} whenever  $1/p-1/q\ge 1/2$ and $(p,q)\neq (2,\infty), (1,2)$. 
 \end{proof}

 \begin{proof}[Proof of Theorem \ref{carl}] {By the argument before Theorem \ref{local} it is easy to see  that \eqref{car} holds only if $1/p-1/q=2/d$. Thus, 
 this and Theorem \ref{local}  shows that \eqref{car} holds if and only if $p$ and $q$ satisfy \eqref{nec-cond} and $1/p-1/q=2/d$, which is equivalent to \eqref{range}. Hence we get  Theorem \ref{carl}.}
\end{proof}

\begin{proof}[Proof of Theorem \ref{dirac}] As before it is sufficient to show that 
$
\normo{\F^{-1}\B(\frac{ (\xi+i v)\wh f(\xi)}{|\xi|^2+2iv\cdot\xi -|v|^2} \B)}_{L^q(\R^2)} \le C\|f\|_{L^p(\R^2)}.
$
By rescaling $\xi\to |v|\xi$  with the condition $1/p-1/q=1/2$  and rotation  this reduces to 
{\[
\normo{\F^{-1}\left(\frac{  R(\xi+ie_2)\wh f(\xi)}{|\xi|^2+2i\xi_2 -1} \right)}_{L^q(\R^2)} \le C\|f\|_{L^p(\R^2)},
\]
with $C$ independent of $R\in \text{SO}(2)$. It is easy to see that $\wt m(\xi)=\frac{|\xi| (1-\chi(\xi)) R(\xi+ie_2)}{|\xi|^2+2i\xi_2 -1}$ satisfies Mihlin's condition uniformly in $R$.} Thus 
$\|\F^{-1} (\wt m\, \widehat g\,)\|_q$ $\lesssim \|g\|_q$ for $1<q<\infty$. Using this and the Hardy-Littlewood-Sobolev inequality we see that, for $p,q$ satisfying $1/p-1/q=1/2$ and $1<p<2$, 
\[\B\| \F^{-1} \left(\frac{(1-\chi(\xi))  R(\xi+ie_2)\wh f(\xi)}{|\xi|^2+2i\xi_2 -1}\right)\B\|_q\lesssim   \B\| \F^{-1} \left(\frac{\wh f(\xi)}{|\xi|}\right)\B\|_q\lesssim \|f\|_p.\]
It remains to show 
\[\B\| \F^{-1} \left(\frac{\chi(\xi)  R(\xi+ie_2)\wh f(\xi)}{|\xi|^2+2i\xi_2 -1}\right)\B\|_q\lesssim \|f\|_p.\]
Since $\chi(\xi)  R(\xi+ie_2)$ is smooth and compactly supported, $\|\F^{-1}\big(\chi(\xi)  R(\xi+ie_2)\wh f(\xi)\big)\|_r\lesssim \|f\|_r$,  $1\le r\le \infty$, independently of $R$. Hence, the estimate follows from Theorem \ref{local}. 
\end{proof}


\section{Non-elliptic cases: Proof of Theorem \ref{non-ell-carleman} }

In this section we are primarily concerned with proving Theorem \ref{non-ell-carleman}.
Sufficiency part follows from the uniform Sobolev estimate \eqref{nonellip-uniform}. Hence it is enough to show that the estimate \eqref{car-non} implies \eqref{nonellip-range}.  

From homogeneity it is easy to see that  the estimate \eqref{car-non} holds only if
\Be\label{homo} \frac1p-\frac1q=\frac 2d\,.\Ee
As before, we note that the estimate \eqref{car-non} is equivalent to 
\begin{equation}\label{mult-non}
\normo{\F^{-1}\B(\frac{\wh f(\xi)}{Q(\xi)+2iv^tM\xi - Q(v)} \B)}_{L^q(\R^d)} \le C\|f\|_{L^p(\R^d)},
\end{equation}
where $Q(\xi)=\xi^t M \xi$ and $M$ is a diagonal matrix with its diagonal entries $-1,\dots, -1,1,\dots, 1$ with $l$-many $-1$s. 
To show that the estimate \eqref{car-non} implies the other conditions of \eqref{nonellip-range}, we  assume that $Q(v)>0$. 
By scaling with the condition \eqref{homo} we may assume $Q(v)=1$. Using { a Lorentz transformation we may also assume  $v=\pm e_{l+1}$.} 
Indeed, by rotations in $(\xi_1, \dots, \xi_l)$-space and $(\xi_{l+1}, \dots, \xi_d)$-space  we may assume $v=(0,\dots, 0,v_l,v_{l+1}, \dots, 0)$ with $-v_l^2+v_{l+1}^2=1$. Then, in the $\xi_l$-$\xi_{l+1}$ plane we perform a change of variables. Let  us set 
\[  \wt M= \begin{pmatrix}  -1 & 0  \\    0 & 1 \end{pmatrix}, \quad  \overline \xi= \begin{pmatrix} \xi_l\\\xi_{l+1}\end{pmatrix}, \quad  R_\theta\overline \xi:= 
  \begin{pmatrix}  \cosh\theta & \sinh\theta  \\    \sinh\theta & \cosh\theta   \end{pmatrix} \overline \xi,\quad \wt{Q}(\overline \xi)=\overline \xi^t \wt M\, \overline \xi.\]
 Since $  Q(v)=\wt Q( (\wt MR_\theta)^t  \overline v )=1$, we can find  $\theta_{\pm}\in \mathbb \R$ such that $(\wt MR_{\theta_\pm})^t  \overline v=(0,\pm1)^t$.  Since $\det(R_{\theta_\pm})=1$, 
 by the change of variables 
 \[(\xi_1, \dots, \xi_{l-1}, \overline \xi^t, \xi_{l+2}, \dots, \xi_{d})\to   (\xi_1, \dots, \xi_{l-1},  (R_{\theta_\pm}\overline \xi)^t, \xi_{l+2}, \dots, \xi_{d}),\]
 we see that the estimate \eqref{mult-non} is equivalent to 
\[
\normo{\F^{-1}\left(\frac{\wh f(\xi)}{Q(\xi) -1 \pm 2i\xi_{l+1}} \right)}_{L^q(\R^d)} \le C\|f\|_{L^p(\R^d)}.
\]
We need to show these estimates hold only if $p,q$ satisfy \eqref{nonellip-range}. 

Let  $ \mathbb  A:=\{ \xi: 1/2\le |\xi|\le 2\}$ and $\beta\in C_0^\infty (\mathbb A)$.  Then, the above {estimate} implies, for any 
$\lambda>0$, 
\[
\normo{\F^{-1}\left(\frac{\beta(\lambda\xi)\wh f(\xi)}{Q(\xi) -1 \pm 2i\xi_{l+1}} \right)}_{L^q(\R^d)} \le C\|f\|_{L^p(\R^d)}.
\]
Together with the condition \eqref{homo}, rescaling $\xi\to \lambda^{-1}\xi$ gives
\[
\normo{\F^{-1}\left(\frac{\beta(\xi)\wh f(\xi)}{Q(\xi) -\lambda^2 \pm 2\lambda i\xi_{l+1}} \right)}_{L^q(\R^d)} \le C\|f\|_{L^p(\R^d)}.
\]
{Considering the difference of the multipliers $(Q(\xi) -\lambda^2 \pm 2\lambda i\xi_{l+1})^{-1}$, we obtain}
\[
\normo{\F^{-1}\left(\left(\frac{2\lambda \xi_{l+1}}{(Q(\xi) -\lambda^2)^2 + 4\lambda^2\xi_{l+1}^2}\right)\,\beta(\xi)\wh f(\xi) \right)}_{L^q(\R^d)} \le C\|f\|_{L^p(\R^d)}.\]
 Note that $\lim_{\lambda\to 0^+} \int \frac{a\lambda }{(t -\lambda^2)^2 + a^2\lambda^2}\phi(t) dt=\pi\phi(0)$ if $a>0$, and take  $f\in \mathcal S(\mathbb R^d)$ of  which  Fourier transform is supported in $\{\xi: \xi_{l+1}\sim 1\}$.  Then,  by letting $\lambda\to 0$ {and Fatou's lemma}, we see that 
\[\B\|\mathcal F^{-1}\B(\chi(\xi)\delta(Q(\xi) )\B)\B\|_{L^q(\R^d)} \le C\]
where $\chi=\widehat f\beta$.
This is possible only if $q>\frac{2(d-1)}{d-2}$. See \cite[pp. 341--343] {JKL}.  Duality gives the condition $p<\frac{2(d-1)}{d}$. This completes the proof.  \hfill \qed


\section{Estimates for the heat equation: Proof of Theorem \ref{carl-heat}} \label{section_heat}

As before, by homogeneity  it is not difficult to see that \eqref{heat} holds only if 
\Be\label{hhomo}\frac 1p-\frac 1q=\frac 2{d+2}.\Ee 
The estimate \eqref{heat} is equivalent to 
\[   \B\|\F^{-1} \left(\frac{\widehat F(\xi,\tau)}{i\tau-\gamma +|\xi|^2+2i\xi\cdot v-|v|^2}\right)\B\|_{q}\le C \|F\|_p. \] 
{By rotation we may set $v=|v|e_d $. By the change of variables $\tau\to \tau-2|v|\xi_d$, which clearly does not {affect} the estimate,   
  the above is  equivalent to }
\[   
\B\|\F^{-1} \left(\frac{\widehat F(\xi,\tau)}{i\tau+|\xi|^2-|v|^2-\gamma}\right)\B\|_{q}\le C \|F\|_p.
\]
Rescaling $( \xi,\tau)\to (  ||v|^2+\gamma|^\frac12 \xi,||v|^2+\gamma|\tau)$ with the condition \eqref{hhomo} reduces the above estimate to
\Be\label{modified} \B\|\F^{-1} \left(\frac{\widehat F(\xi,\tau)}{i\tau+|\xi|^2+\sigma}\right)\B\|_{q}\lesssim \|F\|_p, 
\Ee
where 
\[\sigma=\begin{cases} \,\,\,\,  1&  \text{if} \quad  |v|^2+\gamma<0,   \\
                              \,\,\,\, 0 &  \text{if} \quad  |v|^2+\gamma=0,  \\  
                               -1  &  \text{if} \quad  |v|^2+\gamma>0. \end{cases}  \] 

If $\sigma=0$ or $\sigma=1$, then either the origin is the only singular point of the multiplier or the multiplier has no singularity. Thus,   by direct kernel estimate  it is easy to obtain  \eqref{modified} for $p,q$ satisfying $\frac 1p-\frac 1q=\frac 2{d+2}$, $1<p$, $q<\infty$. This can be handled by making use of analysis in homogeneous spaces. However, we provide an elementary argument. 
Let us set
\[   \|(\xi,\tau)\|=\sqrt{|\xi|^4+\tau^2}. \]  
Then, by scaling and integration by parts we get $|\F^{-1}(\frac{\psi(2^{-j}\|(\xi,\tau)\|)}{i\tau+|\xi|^2+\sigma})| \lesssim 2^{\frac {dj}2 }(1+2^{\frac j2 }|x|+2^j|t|)^{-N}$. Thus  $\| \F^{-1}(\frac{\psi(2^{-j}\|(\xi,\tau)\|)\widehat F(\xi,\tau)}{i\tau+|\xi|^2+\sigma}) \|_1\lesssim 2^{-j} \|F\|_1$ while  $\| \F^{-1}(\frac{\psi(2^{-j}\|(\xi,\tau)\|)\widehat F(\xi,\tau)}{i\tau+|\xi|^2+\sigma})\|_\infty\lesssim 2^{\frac {dj}2} \|F\|_1$. Interpolation between these two estimates {(Bourgain's trick; see \cite{B,CSWaWr, ls0})}  gives the weak type bound 
\[\B\|\F^{-1} \left(\frac{\widehat F(\xi,\tau)}{i\tau+|\xi|^2+\sigma}\right)\B\|_{\frac{d+2}{d}, \infty}\lesssim \|F\|_{1} .\]
Hence, by duality we get the bound from $L^{\frac{d+2}{2},1}\to L^\infty$, and interpolation between these two estimates  gives the desired {bound \eqref{modified}}.  

We finally consider the case $\sigma=-1$. The multiplier now has singularity on the $(d-1)$-sphere $\{(\xi,0)\in\R^d\times\R: |\xi|=1 \}$. Let $\wt\chi$ be a smooth function supported in $\{(\xi,\tau): \|(\xi,\tau) \|\le 4\}$ such that  $\wt\chi=1$ on $\{(\xi,\tau): \|(\xi,\tau) \|\le 2\}$.
By repeating the argument in the above which  deals with the case $\sigma=0,1$, it is easy to see that  for $p,q$ satisfying $\frac 1p-\frac 1q=\frac 2{d+2},$ $1<p,$ $q<\infty$,
\[ \B\|\F^{-1} \left(\frac{(1-\wt\chi(\xi,\tau))\widehat F(\xi,\tau)}{i\tau+|\xi|^2-1}\right)\B\|_{q}\lesssim \|F\|_p. \]
{Thus,  in order to  complete the proof of  the case $\sigma=-1$ it is enough to show the following and this completes the proof of Theorem \ref{carl-heat}. }
\begin{prop}\label{hprop} Let $d\ge1$ and $1<p,q<\infty$. We have the estimate 
\Be\label{hbound}
\B\|\F^{-1} \left(\frac{\wt\chi(\xi,\tau)\widehat F(\xi,\tau)}{i\tau+|\xi|^2-1}\right)\B\|_{q}\lesssim \|F\|_p
\Ee 
if and only if  
\begin{align}
 \label{hcon}
 \frac1p-\frac1q\ge \frac{2}{d+3}, \quad \frac {d+1}p-\frac1q\ge \frac{d+1}{2}, \quad  \frac{d-1}2 \ge  \frac {d+1}q-\frac1p\,.
 \end{align} 
\end{prop}

This can be shown by the exactly same argument which was used for the elliptic case (Theorem \ref{local}). We decompose the multiplier dyadically along $\tau$ and use the restriction estimates to the sphere $\mathbb S^{d-1}$.  Thus, 
we shall be brief.

As before, we set 
\[\mu_\epsilon(\xi,\tau)=\frac{\wt\chi(\xi,\tau)\psi( |\tau|/\epsilon)}{i\tau+|\xi|^2-1}.\] 
Thus, $\frac{\wt\chi(\xi,\tau)}{i\tau+|\xi|^2-1}=\sum_{\epsilon\in \mathbb D} \mu_\epsilon.$  We  make use of Theorem \ref{s-t} and Lemma \ref{rextcy} with $d-2$ replaced by $d-1$, and 
repeating  the same argument in the proof of Lemma \ref{delta-local},  we get the following:

\begin{lem}\label{delta-local-} Let $d\ge 1$.   Suppose  $p, q$ satisfy that  $(\frac1p,\frac1q)$ is contained in the closed triangle $\widetilde{\mathcal T}$ with vertices $\mathfrak S(d+1), {\mathfrak Q(d+1)}, (\frac12,0)$ from which the line segment $[\mathfrak Q(d+1), (\frac12,0))$ is removed. Then we have 
\begin{equation}\label{tri}   \mathfrak M_{p, q} [\mu_\varepsilon]\lesssim \varepsilon^{\frac {d+1}p-\frac1q-\frac {d+1}2}.
\end{equation}
In particular, if $(\frac1p,\frac1q)\in [\mathfrak S(d+1), \mathfrak Q(d+1)\big)$, 
we have
\begin{equation}\label{line-dg1}   \mathfrak M_{p, q} [\mu_\varepsilon]\lesssim \varepsilon^{-1+\frac{d+3}2(\frac1p-\frac1q)}.
\end{equation}
\end{lem}

Making use of the Littlewood-Paley inequality  \eqref{L-P}, by duality and interpolation we get the estimate \eqref{hbound} if \eqref{hcon} holds. For necessity we need to show that 
\[ \mathfrak M_{p, q} [\mu_\varepsilon]\gtrsim 
\max\big(\varepsilon^{\frac {d+1}p-\frac1q-\frac {d+1}2}, \varepsilon^{\frac1p-\frac {d+1}q+\frac {d-1}2}, \varepsilon^{-1+\frac{d+3}2(\frac1p-\frac1q)}\big) {.}
\] 
This can be obtained  by  the same argument in the proof of Lemma \ref{lower}. So we omit the details.

\section{Applications to the unique continuation properties: Proofs of Corollaries  \ref{dirac_unique}, \ref{heat_unique}} 
\label{dirac_sec}

We apply Carleman inequalities (Theorems  \ref{dirac}, \ref{carl-heat}) to prove the unique continuation results (Corollaries  \ref{dirac_unique}, \ref{heat_unique}). The proofs here are similar to those in \cite{KRS}.   Proof of Corollary \ref{heat_unique}  is rather straightforward once one has the estimate \eqref{heat}. So we omit its proof. 
For the proof of Corollary \ref{dirac_unique}   we use the argument  in  \cite{KRS} which is based on reflection principle and the Carleman estimate. For this we use the Kelvin transform of which action on the Laplacian is well known, but the relation between the Kelvin transform and the Dirac operator  does not  seem to be so, although it is likely that  the identities similar to \eqref{k-d} have  been obtained elsewhere. Since we could not find a proper reference for it, we include a proof (see Lemma \ref{id-k}).

\emph{Kelvin transform on the Dirac operator $\mathcal D$ in $\mathbb R^2$. }
We use the explicit representation for the Dirac operator. Let us set 
 \[ D_\pm=-i(\sigma_1\partial_x\pm \sigma_2\partial_y),  \quad  (x,y)\in \mathbb R^2, \] 
 where 
\[\sigma_1=\begin{pmatrix} 0 & 1 \\ 
                                            1&  0\end{pmatrix}, \quad \sigma_2=\begin{pmatrix} 0 & -i \\ 
                                            i&  0\end{pmatrix}.\] 
Then the dirac operator $\mathcal D$ is defined by $\mathcal D= D_+$
(see De Carli-\={O}kaji \cite{DeO}). It is easy to see that $\sigma_1^2=\sigma_2^2=I$,  $\sigma_1\sigma_2=-\sigma_2\sigma_1$ and, thus, $D^2_\pm=-\Delta I$.  Let us set
\[u^\ast(x,y)=(u\circ \Psi)(x,y), \quad   \Psi(x,y)=(X(x,y), Y(x,y))=\left(\frac{x}{x^2+y^2},\frac{y}{x^2+y^2}\right), \]
which is the Kelvin transform. Let us set 
\[  M_\pm(x,y)=\frac{1}{(x^2+y^2)^{2}} \begin{pmatrix}  (ix\pm y)^{2}&0\\ 
                                           0&(ix\mp y)^{2} \end{pmatrix}.\]
\newcommand{\cpsi}{\circ \Psi}
\begin{lem}\label{id-k}  We have the identity 
\Be \label{k-d} D_\pm u^*=  M_\pm ( D_\mp u)^*.
\Ee
\end{lem}

\begin{proof}  Let us note that 
\Be \label{id0} 
X_x'=\frac{y^2-x^2}{(x^2+y^2)^{2}},  \quad   X_y'=Y_x'=\frac{-2xy}{(x^2+y^2)^{2}}, \quad Y_y'=\frac{x^2-y^2}{(x^2+y^2)^{2}}.
\Ee
Then it follows that 
\begin{align}
\label{id}
 i D_\pm u^*&=\sigma_1 \Big((\partial_x u)\circ \Psi\, X_x'+ (\partial_y u)\cpsi\,Y_x' \Big)  \pm \sigma_2\Big((\partial_x u)\circ \Psi\, X_y'+ (\partial_y u)\cpsi\,Y_y') \Big)
\\
\nonumber
 &=\Big(X_x'\sigma_1 \pm  X_y' \sigma_2   \Big) (\partial_x u)\circ \Psi\, +  \Big(Y_x'\sigma_1 \pm   Y_y' \sigma_2   \Big) (\partial_y u)\circ \Psi.
 \end{align} 
 Using \eqref{id0},  we see
\begin{align*} 
X_x'\sigma_1\pm   X_y' \sigma_2
     &=
\begin{pmatrix} 0 & \frac{y^2-x^2}{(x^2+y^2)^{2}}\pm\frac{2ixy}{(x^2+y^2)^{2}}
\\ 
\frac{y^2-x^2}{(x^2+y^2)^{2}} \mp \frac{2ixy}{(x^2+y^2)^{2}} &  0\end{pmatrix}
=
\frac{1}{(x^2+y^2)^{2}} \begin{pmatrix} 0 & (ix\pm y)^{2}\\ 
                                           (ix\mp y)^{2} &  0\end{pmatrix}
                                           \\[6pt]
&=M_\pm (x,y) \sigma_1                                          
                                           \end{align*}
 and                                     
\begin{align*} Y_x'\sigma_1 \pm  Y_y' \sigma_2& =\begin{pmatrix} 0 & \frac{-2xy}{(x^2+y^2)^{2}} \mp i\frac{x^2-y^2}{(x^2+y^2)^{2}}\\ 
                                            \frac{-2xy}{(x^2+y^2)^{2}} \pm i\frac{x^2-y^2}{(x^2+y^2)^{2}}&  0\end{pmatrix}
                                            =
\frac{1}{(x^2+y^2)^{2}} \begin{pmatrix} 0 &\pm i (ix\pm y)^{2}\\ 
                                           \mp i(ix\mp y)^{2} &  0\end{pmatrix}\\[6pt]
                                           &=\mp M_\pm(x,y) \sigma_2 .\end{align*} 
Hence, combining this and \eqref{id} we get the desired identity. 
\end{proof}

Using the Carleman inequality \eqref{carl_dirac} for the Dirac operator $\mathcal D$ and Lemma \ref{id-k} one can prove the following claim which asserts that the zero set of a solution to $|\mathcal Du|\le|Vu|$ is continued (locally) through {a convex curve} in the plane.
\begin{lem}\label{cont_thru_circle}
Let $\Omega\subset\R^2$ be a connected open set containing the unit circle $\mathbb S^1$ and let $V$, $u$ be as in Corollary \ref{dirac_unique}. Suppose that \eqref{diff_dirac} holds in a neighborhood of $\mathbb S^1$. If $u=0$ in the set $\{x\in\R^2: 1<|x|<1+\epsilon_1\}$ for some $\epsilon_1>0$, then there is a positive number $\epsilon_2<1$ such that $u=0$ in the set $\{x: 1-\epsilon_2<|x|<1+\epsilon_1\}$. Conversely, if $u=0$ in $\{x: 1-\epsilon_2<|x|<1 \}$ for some $\epsilon_2$, it follows that $u$ vanishes in $\{x:1-\epsilon_2<|x|<1+\epsilon_1\}$ for some $\epsilon_1$.
\end{lem}
\begin{proof}
First, let us assume that $u=0$ in $\{x: 1<|x|<1+\epsilon_1\}$. By compactness and symmetry it is sufficient to show that $u=0$ in a neighborhood of $-e_2=(0,-1)\in\R^2$. Choose a radial function $\phi\in C_0^\infty(\R^2)$ such that $0\le \phi\le 1$, $\phi(x)=1$ if $|x|\le \epsilon_1/2$, and $\phi(x)=0$ if $|x|\ge \epsilon_1$. Set $\tilde u(x)=\phi(x+e_2)u(x)$. Since $V\in L^2_{loc}(\Omega)$, we can choose $0<r<\epsilon_1/2$ such that 
\[
C\|V\|_{L^2(\{x: |x+e_2|<r, \ -1\le x_2\le -1+\rho \})} \le \frac 12,
\]
where $\rho>0$ denotes the height of the intersection points $\mathbb S^1\cap \{x: |x+e_2|=r\}$, measured from the line $x_2=-1$. Now we apply the Carleman estimate \eqref{carl_dirac} to $\tilde u$ with $v=-\lambda e_2=(0,-\lambda)$, $\lambda>0$. If $1/q=1/p-1/2$\footnote{Trivially, we may assume that $u\in W^{1,p}_{loc}$ for some $1<p<2.$}, then by \eqref{diff_dirac}
\begin{align*}
\|e^{-\la x_2}\tilde u\|_{L^q(\{ x_2\le -1+\rho \})} 
&\le C\|e^{-\la x_2} \mathcal D \tilde u\|_{L^p(\{x_2\le -1+\rho \})} + C\|e^{-\la x_2} \mathcal D\tilde u\|_{L^p(\{ x_2> -1+\rho \})} \\
&\le \frac 12 \|e^{-\la x_2}\tilde u\|_{L^q(\{ x_2\le -1+\rho\})} + C\|e^{-\la x_2} \mathcal D\tilde u\|_{L^p(\{ x_2> -1+\rho \})}.
\end{align*}
From this we have
\[
\|e^{-\la x_2}\tilde u\|_{L^q(\{ x_2\le -1+\rho \})}\le 2C \|e^{-\la x_2} \mathcal D\tilde u\|_{L^p(\{ x_2> -1+\rho \})} \le 2C e^{-\la(\rho-1)}\|\mathcal D\tilde u\|_{L^p}.
\]
If $\lambda$ gets arbitrarily large, this inequality forces $ u$ to vanish on the set $\{x: x_2<-1+\rho,\, |x+e_2|\le r\}$.

In order to prove the other assertion, assume that $u=0$ in $\{x: 1-\epsilon_2<|x|<1\}$. Then $u^*=0$ in $\{x: 1<|x|<1+\epsilon_1\}$ with $\epsilon_1=\frac{\epsilon_2}{1-\epsilon_2}$. Moreover, since $|M_{\pm}(x) v|=|x|^{-2}|v|$, $v\in\R^2$, we have by Lemma \ref{id-k} that 
\[
|\mathcal Du^*(x)|=|x|^{-2}|(D_-u)^* (x)| \le |x|^{-2}|V^*(x) u^*(x)|,
\]
where  $V^\ast=V\circ \Psi$. The last inequality follows from the assumption \eqref{diff_dirac} and the fact that $|\mathcal Du|=|D_-u|$. It is easy to check that $|x|^{-2}V^*(x)\in L^2_{loc}$ in a neighborhood of  $\mathbb S^1$, and by the first case it follows that $u^*=0$ in a neighborhood of the unit circle.
\end{proof}

\begin{proof}[Proof of Corollary \ref{dirac_unique}]
Suppose that $u$ is not identically zero in $\Omega$. Then we can choose a  ball $B_\circ=\{x:|x-x_\circ|<r_\circ\}\subset\Omega$ in which $u=0$ such that the ball is maximal in the sense that $u$ is not identically zero in any larger ball with the same center $x_\circ$. By translation and dilation invariance of the differential inequality \eqref{diff_dirac}, we may assume that $B_\circ=\{x\in\R^2: |x|<1\}$. However, by Lemma \ref{cont_thru_circle}, the zero set of $u$ is in fact larger than $B_\circ$, which is contradiction to the maximality. Therefore $u$ is identically zero in $\Omega$.
\end{proof}


{\bibliographystyle{plain}

\end{document}